\documentclass[11pt,a4paper]{article}

\usepackage[utf8]{inputenc}

\usepackage{amsmath, amsthm, amssymb, amsfonts}

\usepackage{geometry}
\usepackage{parskip} 

\usepackage[numbers,sort&compress]{natbib}

\usepackage{hyperref}
\hypersetup{
    colorlinks=true,
    linkcolor=blue,
    filecolor=magenta,      
    urlcolor=cyan,
    citecolor=red,
}

\newtheorem{theorem}{Theorem}[section]

\newtheorem{definition}[theorem]{Definition}
\theoremstyle{remark}

\theoremstyle{definition}
\newtheorem{example}[theorem]{Example}

\newcommand{\indicator}[1]{\mathbf{1}_{#1}}

\title{\textbf{Eigenvalues of a coupled system of thermostat-type via a Birkhoff--Kellogg type Theorem}}

\author{
    Sajid Ullah\thanks{
       Department of Mathematics and Computer Science, University of Calabria, 
        Ponte P. Bucci 30B, Rende (CS), Italy. 
        \href{mailto:sajid.ullah@unical.it}{sajid.ullah@unical.it}
    }
}
\date{} 
\begin{document}

\maketitle

\begin{abstract}
In this paper, by means of Birkhoff--Kellogg type Theorem in cones we address the existence of eigenvalues and the corresponding eigenvectors to a family of coupled system of thermostat type. The system is characterized by the presence of a real parameter that influences not only the differential equations but also the boundary conditions. Motivated by models of temperature regulation and feedback-controlled systems, we reformulate the original boundary value problems into systems of Hammerstein integral equations. The theoretical results are applied to three different classes of boundary conditions in $t=0$, which are supported by examples.
\end{abstract}

\section{Introduction}
The study of coupled systems of second-order ordinary differential equations (ODEs)
\begin{equation*}
-u_i''(t)=f_i\bigl(t,u_1(t),u_2(t)\bigr),\quad t\in(0,1),\;i=1,2,
\end{equation*}
has become a cornerstone in nonlinear analysis, mostly because it can be used to describe a wide variety of physical and biological processes. An extremely useful and widely-known method of establishing the existence of positive solutions is to rewrite the original system of boundary value problems (BVPs) as a system of perturbed Hammerstein integral equations, which is then studied by topological fixed-point techniques in cones (see \cite{InfantePietramala_NA2009, C.I.P17, I. M. P12, I.P.2014, infante2009eigenvalues, calamai2023affine, agarwal2007positive, I. W2006, Webb2006}). The depth of this area can be seen in the fact that an impressive diversity of boundary conditions has been pursued with success.

A progressive generalization of the boundary structures has been shown in the literature.  Systems with multi-point and four-point coupled BCs are included in fundamental work. For example, Asif and Khan \cite{AsifKhan_JMAA2012} investigate the problem that involved non-homogenous conditions of the form $x(1)=\alpha y(\xi)$ and $y(1)=\beta x(\eta)$.  The Guo–Krasnosel’skiĭ cone-expansion and compression theorem was used to prove their existence.  This classical theorem was also used by Henderson and Luca \cite{HendersonLuca_AMC2012} in the study of a coupled system with multi-point boundary conditions.

The addition of nonlinear and functional boundary conditions marked a major advancement in the generalization process. A flexible framework based on the fixed-point index was developed by Infante and Pietramala \cite{InfantePietramala_NA2009, I.P.2014} to handle systems with fairly general nonlocal and nonlinear BCs, such as those of the type $u_i(0)=H_i[u_j]$ and $u_i(1)=G_i[u_j]$.  Similar techniques were used by Goodrich \cite{Goodrich_AML2015} for systems with the nonlinear BCs $x(0)=H_1[\varphi_1]$ and $y(0)=H_2[\varphi_2]$.  Similarly, Cui and Sun \cite{CuiSun_EJQTDE2012} examined singular superlinear systems with coupled integral BCs $x(1)=\alpha[y] = \int_0^1 y(t)\, dA(t),\ y(1)=\beta[x] = \int_0^1 x(t)\, dB(t)$ and used cone-based techniques to demonstrate their existence. 

These models have gained additional levels of complexity with the passage of time. Infante and Pietramala \cite{InfantePietramala_arXiv2014} generalized their previous work to deal with impulsive BVPs where solutions displaying jumps guided by requirements like $\Delta u_i|_{t=\tau_i}=I_i(u_i(\tau_i))$. Another challenge arises when the nonlinearities involved derivatives. Xu and Zhang \cite{XuZhang_MJM2022} did this by working in a $C^1$ space where they proved existence of an extended Stieltjes integral BCs coupled system by means of the spectral radius of a linear operator associated to it.

Although this rich literature covers a variety of boundary conditions, there are still gaps in the literature pertaining to systems where a parameter $\lambda$ concurrently drives both differential equations and boundary conditions, this structure changes the problem into a nonlinear eigenvalue problem and called the spectrum problem; see Subsection 4.7 of \cite{C.V.Pao} and the references therein. See also \cite{amann1976fixed, infante2009eigenvalues, calamai2023affine} for related work.

Cianciaruso, Infante and Pietramala~\cite[Section 3]{C.I.P17}, by means of classical fixed point theory studied the following problem:
\begin{equation}\label{cip-ex}
    \begin{cases}
     &u''(t) + g_1(t) f_1(t, u(t), v(t)) = 0, \quad \text{a.e. on } [0, 1], \\
    &v''(t) + g_2(t) f_2(t, u(t), v(t)) = 0, \quad \text{a.e. on } [0, 1], \\
    &u'(0) + H_{11}[u,v] = 0,\\
    &u(1) = \beta_1 u(\eta) + H_{12}[u,v], \quad 0 < \eta < 1, \\
    &v'(0) + H_{21}[u,v] = 0,\\ &v(1) = \beta_2 v'(\xi) + H_{22}[u,v], \quad 0 < \xi < 1,
\end{cases}
\end{equation}
In this paper, we address the existence of eigenvalues and the corresponding eigenvectors for a class of parameter-dependent coupled system of thermostat type
similar to~\eqref{cip-ex}. Our motivation for studying such problems lies in their application in physical phenomena, such as modeling the problem of a cooling or heating system controlled by a thermostat. In these heat-flow problems, controllers, reacting to the sensors, are placed in specific points. These are widely studied problems in the context of linear~\cite{I. C2024, C.I2024, C.I2016, C.I.P17,ref7, ref9, Webb2006, I. W2006, ref21, ref23, ref24, ref25, ref26} and nonlinear~\cite{Kalna2002, Kalna2004, Karat2010, Palam2009} controllers. In particular, we study the following system of ODEs
\begin{equation}
\begin{cases}
-u''_1 (t) = \lambda F_1 (t, u_1 (t), u_2 (t)), \quad t \in (0, \, 1),\\
-u''_2 (t) = \lambda F_2 (t, u_2 (t),u_1 (t)), \quad t\in (0,\, 1),
\end{cases}
\label{eq:1}
\end{equation}
subject to three families of coupled functional boundary conditions in which the parameter $\lambda$ also appears:
\begin{enumerate}
    \item[1]- Dirichlet-type conditions in $t=0$:
    \begin{equation}
    \begin{cases}
    u_1(0 ) = \lambda H_1[u_1, u_2 ], \quad \beta_1u'_1(1)+ u_1( \eta_1   ) = \lambda G_1[u_1, u_2],\\
     u_2(0 ) = \lambda H_2[u_1, u_2], \quad \beta_2u'_2(1)+ u_2( \eta_2   ) = \lambda G_2[u_1, u_2],
    \end{cases}
    \label{eq:2}
    \end{equation}
    \item[2]- Neumann-type conditions in $t=0$:
    \begin{equation}
    \begin{cases}
    u'_1(0 ) + \lambda H_1[u_1, u_2]=0, \quad \beta_1u'_1(1)+ u_1( \eta_1   ) = \lambda G_1[u_1, u_2],\\
     u'_2(0 ) +\lambda H_2[u_1, u_2]=0, \quad \beta_2u'_2(1)+ u_2( \eta_2   ) = \lambda G_2[u_1, u_2],
    \end{cases}
    \label{eq:3}
    \end{equation}
     \item[3]- Mixed Neumann and Dirichlet type conditions in $t=0$:
    \begin{equation}
    \begin{cases}
    u'_1(0 ) +\lambda H_1[u_1, u_2]=0, \quad \beta_1u'_1(1)+ u_1( \eta_1   ) = \lambda G_1[u_1, u_2],\\
     u_2(0 ) = \lambda H_2[u_1, u_2],\quad \quad \beta_2u'_2(1)+ u_2( \eta_2   ) = \lambda G_2[u_1, u_2],
    \end{cases}
    \label{eq:4}
    \end{equation}
\end{enumerate}
where $\beta_i>0$, $\beta_i+\eta_i<1$, $\eta_i\in (0, 1)$, $H_i$ and $G_i$ are suitable functionals.
Here we use a Birkhoff-Kellogg type theorem in cones for the existence of eigenvalues for this family of parameter-dependent and coupled BVPs, a setting somewhat similar to the one employed by Infante in~\cite{InfanteBK21} in the context of elliptic systems. In particular, in Section~\ref{sect2} we prove a fairly general eigenvalue and eigenfunction existence result for a class of systems of Hammerstein integral equations (see Theorem~\ref{thm:2}), which covers, as special case the three types of the above mentioned families of systems of BVPs. In Section \ref{sec3} we illustrate explicitly this fact, with the aid of three mathematical examples. Our results are new and complement the previous theory.

\section{Eigenvalues for a system of Hammerstein integral equations}\label{sect2}
We first recall some useful ingredients. Let $(Z,\| \, \|)$ be a real Banach space, a cone $K\subset Z$ is a closed, convex set such that $\alpha K\subset K$ for all $\alpha \geq 0$ and $K\cap (-K)=\{0\}$.
We consider the following sets $$K_R=\{u\in K
: \; \|u\|<R\},\quad \overline{K}_R=\{u\in K
: \; \|u\|\le R\},\quad \partial K_R=\{u\in K
: \; \|u\|=R\},$$ where $R\in (0, +\infty)$.
With these ingredients, we may recall the following Birkhoff-Kellogg type theorem on cones, due to Krasnosel'ski\u{i} and Lady\v{z}enski\u{\i}.
\begin{theorem}[\cite{Krasno, Kra-Lady}]\label{B-K}
    Let $(Z,\| \, \|)$ be a real Banach space, let $\hat{S}:\overline{K}_{R}\to K$ be compact and suppose that $$\inf_{x\in \partial K_{R}}\|\hat{S}x\|>0.$$ Then there exist $\lambda_{0}\in (0,+\infty)$ and $x_{0}\in \partial K_{R}$
such that $x_{0}=\lambda_{0} \hat{S}x_{0}.$
\end{theorem} 

For the application of Birkhoff-Kellogg type theorem in cones, we make the following assumptions on the following system of Hammerstein integral equations. These assumptions are a special case of the ones in~\cite{C.I.P17}.
\begin{equation}\label{eq:inteq}
    u_i(t)=\psi_{i0}H_i[u_1, u_2]+\psi_{i1}G_i[u_1, u_2]+\int_0^1 K_i(t, s)F_i(s, u_1(s),u_2(s))ds, \quad i\in \{1, 2\},
\end{equation}
\begin{itemize}
    \item[(D$_1$)] $F_i : [0, 1]\times [0, +\infty)\times [0, +\infty)\to [0, +\infty)$ hold the Caratheodory-type conditions 
    \begin{itemize}
        \item[(a)] $F_i(\cdot, u_1,u_2)$ is measurable for each fixed $u_1$ and $u_2$ in $[0, +\infty)$,
        \item[(b)] $F_i(t, \cdot,\cdot)$ is continuous for a.e. $t\in [0, 1]$,
        \item[(c)] and for each $R>0$, there exists $\psi_{iR}\in L^{\infty}[0, 1]$ such that $$F_i(t, u_1,u_2)\leq \psi_{iR}(t) \quad \text{for all}\; u_1,\,u_2 \in (0, R)\; \text{and a.e.}\; t\in [0, 1]$$
    \end{itemize}
    \item[(D$_2$)] The kernels $K_i: [0, 1]\times [0, 1]\to [0, +\infty)$ is measurable and for every $\bar{t}\in [0, 1]$, we have $$ \lim_{t\to \bar{t}}|K_i(t, s)-K_i(\bar{t}, s)|=0$$
    \item[(D$_3$)] For every $i=1, 2$ there exist subintervals $[a_i, b_i]\subseteq [0, 1]$, functions $\Phi_i\in L^{\infty}[0, 1]$ and constants $c_i\in (0, 1]$ such that $$ K_i(t, s)\leq \Phi_i(s)\; \text{for}\; t\in [0, 1]\; \text{and a.e.}\; s\in [0, 1] $$ $$ K_i(t, s)\geq c_i\Phi_i(s)\; \text{for}\; t\in [a_i, b_i]\; \text{and a.e.}\; s\in [0, 1] $$  
    \item[(D$_4$)] $\psi_{ij}\in C([0, 1], [0, +\infty))$ and there exists $c_{ij}\in (0, 1]$ such that $$\psi_{ij}(t)\geq c_{ij}\|\psi_{ij}\| \; \text{for every}\,t\in [a_i, b_i]\quad i=1, 2, \; j=0, 1.$$
\end{itemize}

We work in the product space $Y=C[0, 1]\times C[0, 1]$ with the norm $\|(u_1, u_2)\|_Y=\max(\|u_1\|_{\infty},\, \|u_2\|_{\infty})$, here $C[0, 1]$ is Banach space equipped with the infinity norm $\|u_i\|_{\infty}=\sup_{t\in [0, 1]}|u_i(t)|$. We consider a cone $Q$ in $Y$, defined by $$Q=\{ (u_1, u_2)\; \mid \; u_i\in Q_i, \, \text{for}\; i=1, 2\},$$ where $Q_i$ is the cone $$ Q_i=\{x\in C[0, 1]\, :\, x(t)\geq 0 \;\text{for every}\; t\in [0, 1],\; \min_{t\in [a_i, b_i]}x(t)\geq \tilde{c_i}\|x\|_{\infty} \},$$ where $\tilde{c_i}=\min\{c_i, c_{i0}, c_{i1}\}$. We assume that

\begin{itemize}    
    \item[(D$_5$)] $H_i, \,G_i:Q\to [0,+ \infty) $ are compact functionals.
\end{itemize}
Under the assumptions (D$_1$)-(D$_5$), a routine check shows that the integral operator $$S(u_1,u_2):=(S_1(u_1,u_2), S_2(u_1,u_2)),$$ where
\[
\begin{pmatrix}
    S_1(u_1,u_2)\\
    S_2(u_1,u_2)
\end{pmatrix}(t)
:= \begin{pmatrix}
    \psi_{10}(t)H_1[u_1, u_2]+\psi_{11}(t)G_2[u_1, u_2]+\int_0^1 K_1(t, s)F_1(s, u_1(s),u_2(s))ds\\
    \psi_{20}(t)H_2[u_1, u_2]+\psi_{21}(t)G_2[u_1, u_2]+\int_0^1 K_2(t, s)F_2(s, u_1(s),u_2(s))ds
\end{pmatrix},
\]
maps $Q$ into $Q$ and it is compact (see for example Lemma 1 in \cite{I. M. P12}). Note that $$ (u_1, u_2)\in \partial Q_R\implies \begin{cases}
          \|u_i\|_{\infty}\le R \ \text{for some} \ i\in \{1, 2\},\\
          \|u_i\|_{\infty}=R \ \text{for every} \ i\in \{1, 2\}.
      \end{cases}$$ 

For the solvability of the system (\ref{eq:1}) with (\ref{eq:2}), (\ref{eq:3}), or (\ref{eq:4}) we now state and prove the following result.
\begin{theorem}\label{thm:2} Suppose that $R\in (0, +\infty)$ and the following conditions hold for $i=1, 2$
\begin{itemize}
    \item[(1)] There exists $\gamma_{iR}\in C([0, 1],[0, +\infty) $ such that $$ F_i(t, u_1, u_2)\geq \gamma_{iR}(t), \ \text{for all} \; (t, u_1, u_2)\in[a_i,\, b_i]\times  \prod_{k=1}^{2}  [\delta_{ik}\tilde{c}_iR, R],$$ where $\delta_{ik}$ is the classical Kronecker delta function.
    \item[(2)] There exist $\zeta_{iR}^H,\; \zeta_{iR}^G\in [0,+ \infty)$ be such that $$H_i(u_1, u_2)\ge \zeta_{iR}^H \quad\text{and}\quad G_{i}(u_1, u_2)\ge \zeta_{iR}^G,\ \text{for every} \; (u_1, u_2)\in \partial Q_R.$$ 
    \item[(3)] The inequality 
    \begin{equation}\label{C3}
         \sup_{t\in [a_i, b_i]}\Big[\psi_{i0}(t)\zeta_{iR}^H+\psi_{i1}(t)\zeta_{iR}^G + \int_{a_i}^{b_i} K_i(t, s)\gamma_{iR}(s)ds\Big]>0 
    \end{equation}
    holds.
\end{itemize}
  Then there exist $\lambda_R\in (0, +\infty)$ and $(u_{1R}, u_{2R})\in  \partial Q_R$ such that $(u_{1R}, u_{2R})=\lambda_R S(u_{1R}, u_{2R})$.
  \begin{proof}
      As we know, the operator $S$ is compact. We need to prove, $$ \inf_{(u_1, u_2)\in \partial Q_R}\|S(u_1, u_2)\|>0.$$ 
     
Take $(u_1, u_2)\in \partial Q_R$. Let us assume that $$\|u_1\|_{\infty}=R \quad\text{and}\quad \|u_2\|_{\infty}\le R.$$ Then we have
\begin{equation}\label{theq:1}
    \|S(u_1, u_2)\|=\max\{ \|S_1(u_1, u_2)\|_{\infty},  \|S_2(u_1, u_2)\|_{\infty}\}\ge \|S_1(u_1, u_2)\|_{\infty}.
\end{equation}
      Now we have
\begin{multline*}
    \|S_1(u_1, u_2)\|_{\infty}=\sup_{t\in [0, 1]}\left\{\psi_{10}(t)H_1[u_1, u_2]+\psi_{11}(t)G_1[u_1, u_2]+\int_0^1 K_1(t, s)F_1(s, u_1(s),u_2(s))ds\right\} \\
    \ge \sup_{t\in [a_1, b_1]}\left\{\psi_{10}(t)H_1[u_1, u_2]+\psi_{11}(t)G_1[u_1, u_2]+\int_{a_1}^{b_1} K_1(t, s)F_1(s, u_1(s),u_2(s))ds\right\}.
\end{multline*}
Now since $u_1\in Q_1$ and $\|u_1\|_{\infty}=R$, we have 
$$ \tilde{c_1}R\le u_1(t)\le R \quad\text{or}\quad 0\le u_2(t)\le R, \ \text{for every}\; t\in [a_1, b_1].$$
Thus we can use hypothesis (1) and (2). 
\begin{equation}\label{estS1}
      \|S_1(u_1, u_2)\|_{\infty}\ge \sup_{t\in [a_1, b_1]}\left\{\psi_{10}(t)\zeta_{1R}^H+\psi_{11}(t)\zeta_{1R}^G+\int_{a_1}^{b_1} K_1(t, s)\gamma_{1R}(s)ds\right\}.
  \end{equation}
On the other hand, if we assume that  $\|u_1\|_{\infty}\ge R$ and $\|u_2\|_{\infty}= R$, reasoning as above, we obtain the inequality
\begin{equation}\label{estS1}
      \|S_2(u_1, u_2)\|_{\infty}\ge \sup_{t\in [a_2, b_2]}\left\{\psi_{20}(t)\zeta_{2R}^H+\psi_{21}(t)\zeta_{2R}^G+\int_{a_2}^{b_2} K_2(t, s)\gamma_{2R}(s)ds\right\}.
  \end{equation}
In both cases, for every $(u_1, u_2)\in \partial Q_R$ we have, from ~\eqref{theq:1}, that

$$ \|S(u_1, u_2)\|\geq
\min_{i=1,2}  \sup_{t\in [a_i, b_i]}\Bigl\{\psi_{i0}(t)\zeta_{iR}^H+\psi_{i1}(t)\zeta_{iR}^G + \int_{a_i}^{b_i} K_i(t, s)\gamma_{iR}(s)ds\Bigr\}. 
$$
Note that the right-hand side of the above inequality is independent of $(u_1, u_2)$, which implies  
   $$\inf_{(u_1, u_2)\in \partial Q_R}\|S(u_1, u_2)\|> 0, $$
   and the result follows by Theorem~\ref{B-K}.
  \end{proof}
\end{theorem}
\section{Applications to coupled system of BVPs}\label{sec3}
We now proceed to apply Theorem~\ref{thm:2} to the BVPs mentioned in the Introduction.
\subsection{The BVP \eqref{eq:1}-\eqref{eq:2}}
We begin with the BVP
\begin{equation}\label{BVP1}
\begin{cases}
-u''_1 (t) = \lambda F_1 (t, u_1 (t), u_2 (t)),\quad t \in (0, \, 1),\\
-u''_2 (t) = \lambda F_2 (t, u_2 (t),u_1 (t)),\quad t\in (0,\, 1),\\
    u_1(0 ) = \lambda H_1[u_1, u_2 ], \quad \beta_1u'_1(1)+ u_1( \eta_1   ) = \lambda G_1[u_1, u_2],\\
     u_2(0 ) = \lambda H_2[u_1, u_2], \quad \beta_2u'_2(1)+ u_2( \eta_2   ) = \lambda G_2[u_1, u_2].
    \end{cases}
    \end{equation}
The system \eqref{BVP1} can be written in the integral form,
\begin{equation}\label{eq:1-2}
    u_i(t) = \lambda \Big [ \psi_{i,0}(t) H_i[u_1, u_2] + \psi_{i,1}(t) G_i[u_1, u_2] + \int_0^1 K_i(t, s) F_i(s, u_1(s), u_2(s))\, ds\Big], \ i \in \{1, 2\}.
\end{equation}
Hereafter, for simplicity and with slight abuse of notation, we use the same notation for $K_i$, $\Phi_i$, $\psi_{i0}$ and $\psi_{i1}$ for all BVPs. As in ~\cite{I. W2006, C.I2016, C.I2024}, the Green's functions associated with the system are given by:
\begin{align}
    K_i(t,s) &= \frac{t}{\beta_i+\eta_i}  \beta_i +\frac{t}{\beta_i+\eta_i}   \left((\eta_i-s)\indicator{[0,\eta]}(s) \right) - (t-s)\indicator{[0,t]}(s), \label{eq:green_K}
\end{align}
where $\indicator{[a,b]}(x)$ is the indicator function, equal to 1 if $x \in [a,b]$ and 0 otherwise.
With the choice of the subinterval $[a_i,b_i]\subset(0, \beta_i + \eta_i)\subset (0,1)$ for $i \in \{1,2\}$, the hypotheses $(D_2)$ and $(D_3)$  are satisfied. With  $\Phi_i(s)$ is given by
\[
\Phi_i(s) = 
\begin{cases}
s, & \text{if } \beta_i+\eta_i \geq \frac{1}{2}, \\
\left[ \frac{1 - (\beta_i+\eta_i)}{\beta_i+\eta_i} \right] s, & \text{if } \beta_i+\eta_i < \frac{1}{2},
\end{cases}
\]
and the constant $c_i$ is
\[
c_i = 
\begin{cases}
\min \left\{ \frac{a_i\beta_i}{\beta_i+\eta_i}, \frac{\beta_i+\eta_i - b_i}{\beta_i+\eta_i} \right\}, & \text{if } \beta_i+\eta_i \geq \frac{1}{2}, \\[10pt]
\min \left\{ \frac{a_i\beta_i}{1 - (\beta_i+\eta_i)}, \frac{\beta_i+\eta_i - b_i}{1 - (\beta_i+\eta_i)} \right\}, & \text{if } \beta_i+\eta_i < \frac{1}{2}.
\end{cases}
\]
The functions $\psi_{i,0}(t)$, derived in \cite{I. W2006}, are given by:
\begin{align}
\psi_{i,0}(t) &= 1 - \frac{t}{\beta_i+\eta_i}, \label{eq:phi0}
\end{align}
Since $\psi_{i0}$ are decreasing functions on $[0, 1]$, we have $||\psi_{i,0}||_\infty = \psi_{i,0}(0) = 1$. Moreover, for $t \in [a_i,b_i]$, we have: $$\psi_{i,0}(t) \ge \psi_{i,0}(b_i) = 1 - \frac{b_i}{\beta_i + \eta_i}.$$ Thus, condition D$_4$ is satisfied with $c_{i0} = 1 - \frac{b_i}{\beta_i + \eta_i}$. since $b_i< \beta_i + \eta_i$, $c_{i0} \in (0,1)$.

Both $K_i(t, s)$ and $\psi_{i0}$ change sign when $\beta_i + \eta_i<1$, but are positive on the strip $0\le b_i\le t$, $b_i<\beta_i + \eta_i$. For a detailed analysis, we refer the reader to \cite{C.I2016, I. W2006}. 

The functions $\psi_{i,1}(t)$, calculated using the same methodology as \cite{I. W2006}, are: 
\begin{align}
    \psi_{i,1}(t) &= \frac{t}{\beta_i+\eta_i}, \label{eq:phi1}
\end{align}
$\psi_{i,1}(t) = \frac{t}{\beta_i + \eta_i}$ are non-negative and increasing on $[0,1]$. We have $||\psi_{i,1}||_\infty = \psi_{i,1}(1) = \frac{1}{\beta_i + \eta_i}$. For $t \in [a_i,b_i]$, $\psi_{i,1}(t) \ge \psi_{i,1}(a_i) = \frac{a_i}{\beta_i + \eta_i}$. Thus, condition (D$_4$) is satisfied with $c_{i1} = \frac{\psi_{i,1}(a_i)}{||\psi_{i,1}||_\infty} = a_i$. Since $0 < a_i < 1$, $c_{i1} \in (0,1)$.
\begin{definition}
    We say that $\lambda $ is an eigenvalue of the system (\ref{BVP1}), with a corresponding eigenfunction $(u_1,u_2) \in Q$ such that $\|(u_1,u_2)\| > 0$, if the pair $(\lambda, (u_1,u_2) )$ satisfies the system of Hammerstein integral equations (\ref{eq:1-2}).
\end{definition}
Now we can state the following existence result.
\begin{theorem}\label{3.2}
    Let $F_i:[0, 1]\times [0, +\infty)\times [0, +\infty)\to (0, +\infty)$ be continuous and $[a_i, b_i]\subset(0, \beta_i + \eta_i)\subset(0, 1)$. Let $\tilde{c_i}=c_i$ and $R\in (0, +\infty)$, further assume that the conditions (1)-(3) of Theorem \ref{thm:2} hold. Then there exist $\lambda_R$ and $(u_{1R}, u_{2R})\in\partial  Q$ that satisfy the system (\ref{BVP1}).
\end{theorem}
We illustrate the applicability of the previous theorem in a specific example.
\begin{example}
Consider the system
\begin{equation}
\begin{cases}
-u_1''(t) = \lambda\, \dfrac{1}{2} \left( u_1(t) + u_2^3(t) + 2 \right), \\[1.2ex]
-u_2''(t) = \lambda\, \dfrac{1}{2} \left( u_1^2(t) + u_2^2(t) + 1 \right),\\
u_1(0) = \lambda \left( \dfrac{1}{12} u_1\left( \dfrac{1}{3} \right) + \dfrac{1}{12} u_2(1) + \dfrac{1}{3} \right), \\[1.2ex]
u_2(0) = \lambda \left( \dfrac{1}{6} u_1\left( \dfrac{1}{3} \right) + \dfrac{1}{10} u_2(1) + \dfrac{1}{5} \right),\\
    \dfrac{1}{4} u_1'(1) + u_1\left( \dfrac{1}{4} \right) = \lambda \left( \dfrac{1}{2} \left( u_1\left( \dfrac{1}{6} \right) \right)^{1/2} + \dfrac{\sqrt{2}}{20} \left( u_2\left( \dfrac{1}{5} \right) \right)^3 \right),\\
    \dfrac{1}{3}u_2'(1) + u_2\left( \dfrac{1}{4} \right) = \lambda \left( u_1\left( \dfrac{1}{3} \right) + u_2\left( \dfrac{1}{3} \right) \right).
\end{cases}
\label{eq:example_bc2}
\end{equation}

For $R\in (0, +\infty)$, we may take $[a_i, b_i]=\bigg[ \frac{1}{6}, \frac{1}{3}\bigg],$\quad  $\gamma_{1R}(t) = \frac{1}{2}(\tilde{c_1}R+2)$,
     $\gamma_{2R}(t) = \frac{1}{2}((\tilde{c_2}R)^2+1),$ $\tilde{c_1}=\frac{1}{12}$, $\tilde{c_2}=\frac{1}{9}$ $\zeta_{1R}^H = 1/3 > 0$, $\zeta_{2R}^H = 1/5 > 0$, $\zeta_{1R}^G = 0$ and $\zeta_{2R}^G = 0$.
Therefore, Condition 3 of Theorem \ref{thm:2} is satisfied for $i=1$, 
\begin{multline*}
     \sup_{t\in[\frac{1}{6}, \frac{1}{3}]}\Big[(1 - \frac{t}{\beta_1+\eta_1})\frac{1}{3} + \frac{1}{2}(\tilde{c_1}R+2)\int_{\frac{1}{6}}^{\frac{1}{3}} K_1(t, s)ds\Big]\ge\\
      \sup_{t\in[\frac{1}{6}, \frac{1}{3}]}\Big[(\frac{5}{18}) + \frac{1}{24}(R+2)\int_{\frac{1}{6}}^{\frac{1}{3}} c_1 sds\Big]=\\
      \sup_{t\in[\frac{1}{6}, \frac{1}{3}]}\Big[\frac{5}{18} + \frac{1}{3456}(R+2)\Big]>0
\end{multline*}
and $i=2$
\begin{multline*}
    \sup_{t\in[\frac{1}{6}, \frac{1}{3}]}\Big[(1 - \frac{t}{\beta_2+\eta_2})\frac{1}{5} + \frac{1}{2}((\tilde{c_2}R)^2+1)\int_{\frac{1}{6}}^{\frac{1}{3}} K_2(t, s)ds\Big]\ge\\
      \sup_{t\in[\frac{1}{6}, \frac{1}{3}]}\Big[(\frac{4}{35}) + \frac{1}{162}(R^2+1)\int_{\frac{1}{6}}^{\frac{1}{3}} c_2 sds\Big]=\\
      \sup_{t\in[\frac{1}{6}, \frac{1}{3}]}\Big[\frac{4}{35} + \frac{3}{17496}(R^2+1)\Big]>0,
\end{multline*}
which implies that (\ref{C3}) is satisfied for every $R\in (0, +\infty)$.

Thus, we can apply Theorem \ref{thm:2} to obtain our desired solution.
\end{example}
\subsection{The BVP \eqref{eq:1}-\eqref{eq:3}}
The system of BVPs (\ref{eq:1})-(\ref{eq:3}) can also be written in integral form (\ref{eq:1-2}).

As in ~\cite{I. W2006, Webb2006}, the Green's functions associated with the system are given by:
\begin{equation} \label{eq:kernel_K_star}
K_i(t,s) = \beta_i + (\eta_i-s)\indicator{[0,\eta_i]}(s) - (t-s)\indicator{[0,t]}(s).
\end{equation}
With the choice of the subinterval $[a_i,b_i]\subset(0, \beta_i + \eta_i)\subset (0,1)$ for $i \in \{1,2\}$, the hypotheses $(D_2)$ and $(D_3)$  are satisfied. With  $\Phi_i(s)$ is given by
\[ \Phi_{i} (s)= \begin{cases} \beta_i+\eta_i, & \text{if } \;\beta_i+\eta_i \ge  1/2\\ 
1-(\beta_i+\eta_i), & \text{if } \;\beta_i+\eta_i < 1/2 \end{cases} \]
and the constant $c_i$ is
\[
c_{i} = 
\begin{cases}
\displaystyle\frac{\beta_i+\eta_i - b_i}{\beta_i+\eta_i}, & \text{for } \beta_i+\eta_i \geq \frac{1}{2}, \\[10pt]
\displaystyle\frac{\beta_i+\eta_i - b_i}{1 - (\beta_i+\eta_i)}, & \text{for } \beta_i+\eta_i < \frac{1}{2}.
\end{cases}
\]

     From \cite{I. W2006}, we know that 
\begin{equation}\label{eq:psi_iH}
     \psi_{i,0}(t) = \beta_i + \eta_i-t 
\end{equation}
     with $\| \psi_{i,0} \|=\begin{cases}
        \beta_i+\eta_i \quad for \quad \beta_i+\eta_i \ge \frac{1}{2}\\
        1-(\beta_i+\eta_i) \quad for \quad \beta_i+\eta_i< \frac{1}{2}
    \end{cases}.$
    
    $\psi_{i,0}$ are decreasing functions of $t$ for $i\in \{1, 2\}$ and  $\min_{t \in [a_i, b_i]} \psi_{i,0}(t) = \psi_{i,0}(b_i) = \beta_i + \eta_i - b_i$. This minimum is positive if $b_i < \beta_i + \eta_i$.
   
    Therefore, the constant $c_{i0}$ is:
    \[ c_{i0} = \frac{\min_{t \in [a_i, b_i]} \psi_{i,0}(t)}{\| \psi_{i,0} \|_\infty} = \frac{\beta_i + \eta_i - b_i}{\beta_i + \eta_i}. \]
In this case we calculated the functions $\psi_{i,1}$, as 
\[
\psi_{i,1}:= 1,
\]
 which are non-negative on $[0,1]$. We have $\| \psi_{i,1} \|_\infty = 1$.
    On the subinterval $[a_i, b_i]$, $\min_{t \in [a_i, b_i]} \psi_{i,1}(t) = 1$.
    Thus, the constant $c_{i1}$ is:
    \[ c_{i1} = \frac{\min_{t \in [a_i, b_i]} \psi_{i,1}(t)}{\| \psi_{i,1} \|_\infty} = \frac{1}{1} = 1. \]
Now for the solvability of the system (\ref{eq:1})-(\ref{eq:3}), we can state a result similar to Theorem \ref{3.2}.
\begin{example}
Consider the system
    \begin{equation}
\begin{cases}
-u_1''(t) = \lambda\, \left( (u_1(t))^2 + sin^2(u_2(t)) + 1 \right), \\[1.2ex]
-u_2''(t) = \lambda\,\left( e^{u_1(t)} + u_2^3(t) + 1 \right),\\
u'_1(0) +\lambda \left( \dfrac{1}{10} u_1\left(1 \right) + \dfrac{1}{10} u_2(1) + \dfrac{1}{5} \right)=0, \\[1.2ex]
u'_2(0) + \lambda \left( \dfrac{1}{10} u_1\left( \dfrac{1}{2} \right) + \dfrac{1}{20} u_2(1) + \dfrac{1}{10} \right)=0,\\
    \dfrac{1}{4} u_1'(1) + u_1\left( \dfrac{1}{4} \right) = \lambda \left( \dfrac{1}{4} \Bigl( \sqrt{u_1\Bigl( \dfrac{1}{4} \Bigr)} \Bigr) + \dfrac{1}{8} \left( u_2(1)^2 \right)+\dfrac{1}{5} \right),\\
    \dfrac{1}{3} u_1'(1) + u_1\left( \dfrac{1}{4} \right) = \lambda \left( \dfrac{1}{6}u_1\left( \dfrac{1}{3} \right) + \dfrac{1}{6}u_2\left( \dfrac{1}{4} \right) +\dfrac{1}{6}\right).
\end{cases}
\label{eq:example_bc4}
\end{equation}
For $R\in (0, +\infty)$, we may take $[a_i, b_i]=\bigg[ \frac{1}{6}, \frac{1}{3}\bigg]$, $\gamma_{1R}(t) = (\tilde{c_1}R)^2 +1)$,
     $\gamma_{2R}(t) = (\tilde{c_2}R)^2 +2),$ $\tilde{c_1}=\frac{1}{3}$, $\tilde{c_2}=\frac{3}{7}$, $\zeta_{1R}^H =\frac{1}{5}$, $\zeta_{2R}^H =\frac{1}{10}$, $\zeta_{1R}^G =\frac{1}{5}$ and $\zeta_{2R}^G =\frac{1}{6}$.
     Therefore, Condition 3 of Theorem \ref{thm:2} is satisfied for $i=1$, 
\begin{multline*}
    \sup_{t\in [\frac{1}{6}, \frac{1}{3}]}\Big[(\beta_1 + \eta_1-t)\zeta_{1R}^H+\zeta_{1R}^G + \int_{\frac{1}{6}}^{\frac{1}{3}} K_1(t, s)\gamma_{1R}(s)ds\Big]\ge\\
    \sup_{t\in [\frac{1}{6}, \frac{1}{3}]}\Big[(\beta_1 + \eta_1-b_1)\frac{1}{5}+\frac{1}{5} + \frac{(\tilde{c_1}R)^2 +1)}{2}\int_{\frac{1}{6}}^{\frac{1}{3}} c_1 ds\Big]=\\
    \sup_{t\in [\frac{1}{6}, \frac{1}{3}]}\Big[\frac{7}{30} + \frac{(R)^2 +9)}{324}\Big]>0
\end{multline*}
and for $i=2$
\begin{multline*}
    \sup_{t\in [\frac{1}{6}, \frac{1}{3}]}\Big[(\beta_2 + \eta_2-t)\zeta_{2R}^H+\zeta_{2R}^G + \int_{\frac{1}{6}}^{\frac{1}{3}} K_2(t, s)\gamma_{2R}(s)ds\Big]\ge\\
    \sup_{t\in [\frac{1}{6}, \frac{1}{3}]}\Big[(\beta_2 + \eta_2-b_2)\frac{1}{10}+\frac{1}{6} + \frac{(\tilde{c_2}R)^2 +2)}{2}\int_{\frac{1}{6}}^{\frac{1}{3}} c_2 ds\Big]\ge\\
    \sup_{t\in [\frac{1}{6}, \frac{1}{3}]}\Big[\frac{37}{210} + \frac{9R^2+98}{1372}\Big]>0,
\end{multline*}
which implies that (\ref{C3}) is satisfied for every $R\in (0, +\infty)$. Thus, we can apply Theorem \ref{thm:2} to obtain our desired solution.
\end{example}
\subsection{The BVP \eqref{eq:1}-\eqref{eq:4}}
The system (\ref{eq:1}), together with the set of nonlinear and nonlocal boundary conditions (\ref{eq:4}) can also be written in the form of system of integral equations (\ref{eq:1-2}). 
Where functions $\psi_{i,0}(t), \,\psi_{i,1}(t)$ and $K_i(t,s)$ are given by:
\begin{align*}
\psi_{1,0}(t) &= \beta_1 + \eta_1-t \\
\psi_{1,1}(t) &= 1 \\
\psi_{2,0}(t) &= 1 - \frac{t}{\beta_2+\eta_2}, \label{eq:phi0} \\
\psi_{2,1}(t) &= \frac{t}{\beta_2+\eta_2},\\
K_1(t,s) &= \beta_1 + (\eta_1-s)\indicator{[0,\eta_1]}(s) - (t-s)\indicator{[0,t]}(s),\\
K_2(t,s) &= \frac{t}{\beta_2+\eta_2}  \beta_2 +\frac{t}{\beta_2+\eta_2} (\eta_2-s)\indicator{[0,\eta_2]}(s)  - (t-s)\indicator{[0,t]}(s). 
\end{align*}
The rest is similar to the previous discussion.
\begin{example} Consider the BVP
    \begin{equation}
\begin{cases}
-u_1''(t) = \lambda\, \left( (u_1(t))^2 + \sin^2(u_2(t)) + 1 \right), \\[1.2ex]
-u_2''(t) = \lambda\, \dfrac{1}{2} \left( u_1^2(t) + u_2^2(t) + 1 \right),\\
u'_1(0) +\lambda \left( \dfrac{1}{10} u_1\left(1 \right) + \dfrac{1}{10} u_2(1) + \dfrac{1}{5} \right)=0, \\[1.2ex]
u_2(0) = \lambda \left( \dfrac{1}{6} u_1\left( \dfrac{1}{3} \right) + \dfrac{1}{10} u_2(1) + \dfrac{1}{5} \right),\\
    \dfrac{1}{4} u_1'(1) + u_1\left( \dfrac{1}{4} \right) = \lambda \left( \dfrac{1}{4} \left( \sqrt{u_1\left( \dfrac{1}{4} \right)} \right) + \dfrac{1}{8} \left( u_2(1)^2 \right)+\dfrac{1}{5} \right),\\
    \dfrac{1}{3}u_2'(1) + u_2\left( \dfrac{1}{4} \right) = \lambda \left( u_1\left( \dfrac{1}{3} \right) + u_2\left( \dfrac{1}{3} \right) \right).
\end{cases}
\end{equation}
For $R\in (0, +\infty)$, we can take $[a_i, b_i]=\bigg[ \frac{1}{6}, \frac{1}{3}\bigg]$, $\gamma_{1R}(t) = (\tilde{c_1}R)^2 +1)$,
     $\gamma_{2R}(t) =\frac{1}{2}(\tilde{c_2}R)^2 +1) ,$ $\tilde{c_1}=\frac{1}{3}$, $\tilde{c_2}=\frac{1}{9}$, $\zeta_{1R}^H = 1/5 > 0$, $\zeta_{2R}^H = 1/5$, $\zeta_{1R}^G = 1/5$ and $\zeta_{2R}^G = 0$.
 Therefore, Condition 3 of Theorem \ref{thm:2} is satisfied for $i=1$, 
\begin{multline*}
    \sup_{t\in [\frac{1}{6}, \frac{1}{3}]}\Big[(\beta_1 + \eta_1-t)\zeta_{1R}^H+\zeta_{1R}^G + \int_{\frac{1}{6}}^{\frac{1}{3}} K_1(t, s)\gamma_{1R}(s)ds\Big]\ge\\
    \sup_{t\in [\frac{1}{6}, \frac{1}{3}]}\Big[(\beta_1 + \eta_1-b_1)\frac{1}{5}+\frac{1}{5} + \frac{(\tilde{c_1}R)^2 +1)}{2}\int_{\frac{1}{6}}^{\frac{1}{3}} c_1 ds\Big]=\\
    \sup_{t\in [\frac{1}{6}, \frac{1}{3}]}\Big[\frac{7}{30} + \frac{(R)^2 +9)}{324}\Big]>0
\end{multline*}
and $i=2$
\begin{multline*}
    \sup_{t\in[\frac{1}{6}, \frac{1}{3}]}\Big[(1 - \frac{t}{\beta_2+\eta_2})\frac{1}{5} + \frac{1}{2}((\tilde{c_2}R)^2+1)\int_{\frac{1}{6}}^{\frac{1}{3}} K_2(t, s)ds\Big]\ge\\
      \sup_{t\in[\frac{1}{6}, \frac{1}{3}]}\Big[(\frac{4}{35}) + \frac{1}{162}(R^2+1)\int_{\frac{1}{6}}^{\frac{1}{3}} c_2 sds\Big]=\\
      \sup_{t\in[\frac{1}{6}, \frac{1}{3}]}\Big[\frac{4}{35} + \frac{3}{17496}(R^2+1)\Big]>0,
\end{multline*}

A similar result as the ones above also holds in this example.

\textbf{Acknowledgement} The author is deeply grateful to Professor Gennaro Infante for his insightful comments and constructive suggestions, which have significantly improved the quality of this work.
\end{example}

\end{document}